\documentclass[11pt]{article}
\usepackage{latexsym}
        \usepackage{amscd}
        \usepackage{amssymb}
\usepackage{oldgerm}
        \usepackage{amsmath}
        \usepackage{amsthm}
        \usepackage{enumerate}
\def\co{\colon\thinspace}

\def\e{\epsilon}
\def\a{\alpha}

\newtheorem{thm}{Theorem}[section]

\newtheorem{lem}[thm]{Lemma}

\newtheorem{defn}[thm]{Definition}

\newtheorem{Example}[thm]{Example}
\newenvironment{ex}{\begin{Example}\rm}{\end{Example}}
\newtheorem{Counterexample}[thm]{Counterexample}

\newtheorem{remark}[thm]{Remark}
\newenvironment{rmk}{\begin{remark}\rm}{\end{remark}}
\newtheorem{Fact}[thm]{Fact}

\newtheorem{Nothing}[thm]{$\!\!\!$}

\newcommand{\be}{\begin{equation}}
\newcommand{\ee}{ \end{equation}} 
\newcommand{\ba}{\begin{eqnarray}}
\newcommand{\ea}{\end{eqnarray}}
\newcommand{\ban}{\begin{eqnarray*}}
\newcommand{\ean}{\end{eqnarray*}}

\newcommand{\scal}{\mathrm{scal}}
\newcommand{\Ric}{\mathrm{Ric}}
\newcommand{\Iso}{\mbox{Iso}}

\begin{document}
\abovedisplayskip=6pt plus3pt minus3pt
\belowdisplayskip=6pt plus3pt minus3pt
\title{\bf Nonnegative curvature, symmetry and fundamental group
\footnotetext{\it 2000 Mathematics Subject classification.\rm\
Primary 53C21, 57S15.
Key words: nonnegative curvature, isometry group, splitting,
equivariant connected sum.}\rm}
\author{Igor Belegradek
\thanks {Partially supported by NSF Grant \# DMS-0203979.}}
\date{}
\maketitle
\begin{abstract}
We prove a result on equivariant deformations of 
flat bundles, and as a corollary, we obtain  
two ``splitting in a finite cover'' theorems
for isometric group actions on Riemannian manifolds 
with infinite fundamental groups, where the manifolds are either 
compact of $\Ric\ge 0$, or complete of $\sec\ge 0$.
\end{abstract}

\section{Introduction}
\label{sec: intro}
 J.~Cheeger and D.~Gromoll~\cite{CG2} showed,
as an application of their Splitting Theorem, that
if $M$ is a compact manifold with $\Ric\ge 0$ and 
infinite $\pi_1(M)$, then $M$ has a structure of a flat
(orbifold) bundle whose generic 
fiber is a simply-connected manifold $C$ with $\Ric\ge 0$,
the base is a flat orbifold,
and the holonomy group lies in $\Iso(C)$.
In particular, $M$ has a finite cover 
diffeomorphic to the Riemannian product 
of a flat torus $T$ and 
a simply-connected manifold $C$ with $\Ric\ge 0$.
The cover $C\times T\to M$ is Riemannian
(i.e. a local isometry) precisely when the holonomy group of the
flat bundle is finite. 
For example, if $M$ is Ricci flat, 
the cover is Riemannian, because in this case $C$ is
Ricci flat and so $\Iso(C)$ is finite by Bochner's theorem.
Now if the cover $C\times T\to M$ is Riemannian, then
one can understand the group $\Iso(M)$ by relating it to 
$\Iso(C\times T)\cong\Iso(C)\times\Iso(T)$;
this turns out to be true even for non-Riemannian covers, as follows.

\begin{thm}\label{ric-deform}
Let $(M,g)$ be a compact Riemannian manifold of $\Ric\ge 0$
with $\pi_1(M)$ infinite,
and let $H$ be a compact subgroup of $\Iso(M,g)$. Then there exists
a smooth deformation of metrics $g_t$ on $M$ such that $g_0=g$,
$H\le\Iso(M,g_t)$,
the universal Riemannian covers of $(M,g_t)$ are all isometric, 
and a finitely-sheeted Riemannian cover of $(M, g_1)$
is isometric to the Riemannian product of a flat torus $T$ and
a simply-connected manifold $C$ of $\Ric\ge 0$. 
\end{thm}

Theorem~\ref{ric-deform} generalizes a result of 
B.~Wilking~\cite{Wil} who proved it for trivial $H$. 
Theorem~\ref{ric-deform} follows from a general result
on equivariant deformations of flat bundles
(see Section~\ref{sec: eq split} for details).
Essentially, we show that
any flat bundle of this kind can be equivariantly
deformed, through flat bundles, to a flat bundle with finite holonomy.
This also gives a version of Theorem~\ref{ric-deform} 
for isometric actions on 
complete manifolds of $\sec\ge 0$. Since the isometries of $T$
are well-understood,
Theorem~\ref{ric-deform} roughly means that studying
the isometries  $M$ of
can be reduced to studying the isometries of $C$. We now list some
sample corollaries of Theorem~\ref{ric-deform}.
\begin{itemize}
\item
If $H$ is connected and simply-connected, then the $H$-action on $(M, g_1)$
lifts uniquely to an isometric $H$-action on $C\times T$
such that the $H$-action on the $T$-factor is trivial.
\item
If $H$ has a fixed point $x$ and the induced $H$-action on
$\pi_1(M,x)$ is trivial, then the $H$-action on $M$
lifts to an isometric $H$-action 
on $C\times T$ such that the $H$-action on the $T$-factor is trivial.
\item
If $H$ is connected, then the $H$-action on $M$ lifts to an isometric
action of a finite cover $\bar H$ of $H$ on $C\times T$. 
In this case the $\bar H$-action on the $T$-factor is also 
well-understood 
(see Appendix~\ref{actions on tori}).
\end{itemize}

Given a compact Lie group $G$ acting smoothly and 
effectively on a closed manifold $M$, it is natural to ask whether
$M$ admits a $G$-invariant metric satisfying $\Ric\ge 0$ or $\Ric >0$.
It is well-known that if $G$ acts transitively, then it preserves a 
metric of $\Ric\ge 0$, which has $\Ric>0$ if and only if 
$\pi_1(M)$ is finite. 
It was shown in~\cite{GZ} that the same is true for cohomogeneity 
one actions.
By~\cite{GPT}, $M$ admits a $G$-invariant metric of $\Ric> 0$
if the action is free, $G$ is connected, $M/G$ admits a metric of $\Ric>0$, 
and $\pi_1(M)$ is finite. 
However, for general actions little is known.
Theorem~\ref{ric-deform} gives rise to many examples of 
group actions on compact manifolds that are not isometric in 
any metric of $\Ric\ge 0$. For instance, we prove the following.

\begin{thm}\label{ric-scal-eq-sum} 
Let $M$ be a compact Riemannian 
manifold with $\pi_1(M)$ infinite,
and let $H$ be a closed subgroup of $\Iso(M)$.
Suppose the set of the fixed points of $H$ has a component $F$
of dimension $\ge 3$ and for some $x\in F$, 
the $H$-action on $\pi_1(M, x)$ is trivial.
Then there are infinitely many smooth $H$-actions on $M$
that are non-isometric in any metric of $\Ric\ge 0$, and
such that if $\scal(M)>0$ and $\dim(F)\ge 5$, then the actions are
isometric in some metrics of $\scal>0$.
\end{thm}

To date, there seem to be no examples of smooth $G$-action 
on a compact manifold $M$ with {\it finite} $\pi_1(M)$ 
such that the $G$-action preserves a metric of $\scal\ge 0$
and preserves no metric of $\Ric\ge 0$. 
This problem is well-known for trivial $G$, and one wonders
whether it could become easier for nontrivial $G$. 
If the assumption $\Ric\ge 0$ is replaced by $\sec\ge 0$,
we construct such examples as follows.

\begin{thm}\label{sec-scal-eq-sum} 
Let $M$ be a complete Riemannian manifold
and let $H$ be a compact subgroup of $\Iso(M)$.
Suppose the set of the fixed points of $H$ has 
a component $F$ of dimension $\ge 3$, and for some $x\in F$, 
the $H$-action on $\pi_1(M, x)$ is trivial.
Then there are infinitely many smooth $H$-actions on $M$
that are non-isometric in any complete metric of $\sec\ge 0$, and
such that if $\scal(M)>0$ and $\dim(F)\ge 5$, then the actions are
isometric in some complete metrics of $\scal>0$.
\end{thm}

Note that the question whether a given $G$-action on $M$ 
preserves a metric of $\scal> 0$ is also wide open, except for
the case of free actions: by O'Neill's formula, 
if $G$ acts freely, and either $M/G$ admits 
a metric of $\scal>0$, or $G$ has a biinvariant metric 
of $\scal>0$, then $M$ has a $G$-invariant metric
of $\scal>0$.
See also~\cite{Lot} for obstructions in the case of 
semifree circle actions.

The $H$-actions constructed in 
Theorems~\ref{ric-scal-eq-sum}--\ref{sec-scal-eq-sum} are obtained
from the original action by taking an equivariant connected sum
at a fixed point with certain $H$-actions on the standard sphere.
There are other ways to produce $H$-actions with the same properties
as in Theorems~\ref{ric-scal-eq-sum}--\ref{sec-scal-eq-sum}, 
but for the sake of brevity we do not consider them here, and leave
it to the interested reader to pursue.

This paper grew out of the author's joint work~\cite{BK1} with Vitali 
Kapovitch. I am grateful to
Slawomir Kwasik for Remark~\ref{h-cobordism}, and to Reinhard Schultz
for the idea of using equivariant connected sums in
Theorems~\ref{ric-scal-eq-sum}--\ref{sec-scal-eq-sum}.

\section{Equivariant splitting in a finite cover}
\label{sec: eq split}

Cheeger and Gromoll proved~\cite{CG2} that a closed Riemannian 
manifold $(M,g)$ of $\Ric\ge 0$ has a finite cover diffeomorphic 
to the product of flat torus and a simply-connected manifold of 
$\Ric\ge 0$. In general, the pullback of $g$ to the finite cover
is not a product metric, e.g. the metric on $S^2\times S^1$
obtained by pushing down the standard product metric on
$S^2\times \mathbb R$ to $(S^2\times\mathbb R)/\mathbb Z$
with the $\mathbb Z$-action given by $n\cdot (v,t)=(e^{in}v,t+n)$
has $\sec\ge 0$ and is not a product metric in any finite 
cover, as was noted in~\cite{CG1}.

However, Wilking showed in~\cite{Wil} 
the metric $g$ can be deformed, through metrics of $\Ric\ge 0$, 
to a metric $g_1$ 
such that a finitely-sheeted Riemannian cover of $(M, g_1)$
is isometric to the Riemannian product of a flat torus and
a simply-connected manifold with $\Ric\ge 0$.
A similar conclusion holds in some other situations
e.g. when $M$ is a complete manifolds of $\sec\ge 0$~\cite{Wil, BK1}.

In this section we show, by a modification of Wilking's argument, 
that one can choose the deformation from $g$ to $g_1$ to be 
equivariant with respect to the isometry group of $(M,g)$.
Our exposition is self-contained yet for the reader's convenience we
mostly follow Wilking's notations. Like in~\cite{Wil}, we actually 
prove a much more general result with no curvature assumptions,
stated in the following definition.

\begin{defn}\label{def of star}
If $M$ be a complete Riemannian manifold 
and $H$ is a compact subgroup of the isometry group of $M$, then
the pair $(M,H)$ is said to satisfy $\mathrm{(\ast)}$
if there exists a regular Riemannian covering $\tilde M\to M$
with the infinite group of covering transformation $\Pi$
such that \newline
(i) any element of $H$ lifts to 
an isometry of $\tilde M$, and\newline
(ii) if $N\times\mathbb R^n$
is the de Rham decomposition of \rm $\tilde M$~\cite{EH}\it, 
where $N$ is the product of the non-Euclidean factors, then
the $\Pi$-action on $N$, which is the composition of 
the deck-transformation action 
$\rho\co\Pi\to\Iso(\tilde M)\cong\Iso(N)\times\Iso(\mathbb R^n)$
with the projection onto the $\Iso(N)$-factor, has a 
precompact image.
\end{defn}

If $(M,H)$ satisfies $\mathrm{(\ast )}$, we use the following 
notations.  
We denote the metrics on $M$, $\tilde M$ by $g$, $\tilde g$,
respectively. 
We write the isomorphism 
$\Iso (N\times\mathbb R^n)\cong\Iso (N)\times\Iso (\mathbb R^n)$
as $k\to (k_n, k_r)$ for 
$k\in \Iso (N\times\mathbb R^n)$, and given a subgroup $K$ of
$\Iso (N\times\mathbb R^n)$, we denote by $K_n$, $K_r$
the projections of $K$ into $\Iso (N)$, $\Iso (\mathbb R^n)$,
respectively. Thus, condition (ii) of Definition~\ref{def of star} 
means that $\rho(\Pi)_n$ is precompact in $\Iso(N)$.
Let $\tilde H$ be the set of all the lifts of elements of $H$. 
By~\cite[I.Theorem 9.3]{Bre}, $\tilde H$ 
is a Lie subgroup of 
$\Iso(\tilde M)$ and the map $\tilde H\to H$ sending $\tilde h$ to $h$,
where $\tilde h$ is a lift of $h$,
is a surjection of Lie groups with kernel $\Pi$.

\begin{rmk}
Since $\Pi$ is infinite and discrete, (ii) implies that $n>0$.
Also (i) is a purely topological condition on the $H$-action on $M$,
which for example is true automatically if $\tilde M$ is the universal 
cover of $M$.
\end{rmk}

\begin{ex} \label{ex of cond star}
$(M,H)$ satisfies $\mathrm{(\ast )}$ if 
$\tilde M$ is the universal cover of $M$, $\pi_1(M)$ 
is infinite, $H$ is a compact subgroup of $\Iso(M)$, and
$M$ satisfies one of the conditions below.\newline
(1) $M$ is a compact of $\Ric\ge 0$. Here $\rho(\Pi)_n$ is precompact
since $N$ is compact~\cite{CG2};\newline 
(2) $M$ is a complete of $sec\ge 0$. Here $\rho(\Pi)_n$ is precompact
since $\Iso(N)$ is compact~\cite{CG1};\newline
(3) $M$ is complete of $\Ric\ge 0$ and $M$ 
is the normal bundle of a compact submanifold $S$ such that
either $S$ is totally convex, or there exists a distance
nonincreasing retraction $N\to S$. 
Here $\Iso(N)$ may be noncompact, yet $\rho(\Pi)_n$ is precompact
as is shown in~\cite[Proposition 2.2]{BK1}, where a more general
result is proved. 
\end{ex}
\begin{rmk}
By~\cite{CG1, Per}, if $M$ satisfies~\ref{ex of cond star}(2),
then it satisfies~\ref{ex of cond star}(3).
\end{rmk}

\begin{rmk}
Topologically, if $(M,H)$ satisfies $\mathrm{(\ast)}$, then
$M$ is a flat (orbifold) bundle with holonomy group
$\rho(\Pi)_n$, 
the generic fiber $N$, and the base $\mathbb R^n/\rho(\Pi)_r$.
The $H$-action on $M$ takes fibers to fibers. 
Theorem~\ref{thm: deform metric}, the main result of this section,
gives an $H$-equivariant deformation of such flat
bundles that starts with the given bundle and ends in a bundle with 
a finite holonomy group.
\end{rmk}

\begin{rmk}
The main difference between the proof in~\cite{Wil}
and the arguments of this section is that
Wilking deforms the
deck-transformation group $\Pi$ in $\Iso(\tilde M)$, while 
we deform the whole group $\tilde H$, which surjects onto $H$
with kernel $\Pi$. The compactness
of $H$ comes into the proof on several occasions. 
\end{rmk}

\begin{lem}\label{lem: stabilize}
If $(M,H)$ satisfies $\mathrm{(\ast )}$, then
there is an $d$-dimensional affine subspace $V^d\subset\mathbb R^n$
with $d\in [1,n]$, such that $\tilde H_r$ stabilizes $V^d$ and 
$\Pi$ acts on $V^d$ 
with finite kernel and the image
is a discrete cocompact subgroup of $\Iso(V^d)\cong\Iso (\mathbb R^d)$.
In particular, $\Pi$ has a finite index normal subgroup isomorphic to
$\mathbb Z^d$.
\end{lem}
\begin{proof}
Since $\rho(\Pi)_n$ is precompact in $\Iso(N)$,
the homomorphism $\Pi\to \rho(\Pi)_r$,
which is the composition of $\rho$ with the projection to
$\Iso(\mathbb R^n)$,
has finite kernel and $\rho(\Pi)_r$
is a discrete subgroup of $\Iso(\mathbb R^n)$. 
In particular, $\Pi$ has a normal finite index subgroup $\Pi_1$
isomorphic to $\mathbb Z^d$ for some $d\in [1,n]$,
and furthermore, 
we arrange the subgroup to be invariant under conjugation
in $\tilde H$. (The conjugation by elements of $\tilde H$ 
defines a homomorphism 
$\tilde H\to\mathrm{Aut}(\Pi)$, and since $\Pi$ is finitely 
generated, any finite index subgroup of $\Pi$
contains a finite index subgroup that is invariant under
all the endomorphisms of $\Pi$~\cite[Exercise 15.2.3]{KP}).

The group $\Pi_1$ acts on $\mathbb R^n$ by freely, hence by
the Soul Theorem, which in this case can be found 
in~\cite[Theorem 3.3.3]{Wol}, $\Pi_1$
acts cocompactly on an affine $d$-dimensional subspace 
$W^d\subset \mathbb R^n$. 
Since $\Pi_1\unlhd\tilde H$, for each
$\gamma\in\Pi_1$ and $\tilde h\in\tilde H$, the transformation
$\gamma_r$ stabilizes $\tilde h_r(W^d)$.
Hence $\tilde h_r(W^d)$ is an affine subspace 
parallel to $W^d$.
%
%
Therefore, the $\tilde H$-action on $\mathbb R^n$ descends to an isometric 
$\tilde H$-action on $\mathbb R^n/W^d$, which is
isometric to $\mathbb R^{n-d}$, and $\Pi_1$
lies in the kernel of the action. Now $\tilde H/\Pi_1$
is compact and any isometric action of a compact Lie group on 
a Euclidean space has a fixed point. Thus, $\tilde H$
stabilizes an affine subspace $V^d$ of $\mathbb R^n$ parallel to
$W^d$. 
Note that $\Pi$ acts on $V^d$ with finite kernel and the image
is a discrete cocompact subgroup of $\Iso(V^d)\cong\Iso (\mathbb R^d)$.
\end{proof}

The stabilizer of $V^d$ in $\Iso(\mathbb R^n)$
can be identified to $O(n-d)\times\Iso (\mathbb R^d)$,
so by Lemma~\ref{lem: stabilize}, 
we can view $\tilde H_r$ as a subgroup of 
$O(n-d)\times\Iso (\mathbb R^d)$. 
Since $H$ is compact and $\rho(\Pi)_n$ is precompact,
$\tilde H_n$ is also precompact. Then the closure of 
$\tilde H_n\times O(n-d)$ in $\Iso(\tilde N)\times O(n-d)$
is a compact Lie subgroup, which we denote by $G$.
Thus, $\tilde H$ can be viewed
as a subgroup of $G\times\Iso (\mathbb R^d)$, and
we denote the inclusion by 
$(\psi, \e)\co\tilde H\to
G\times\Iso (\mathbb R^d)$.

\begin{thm}\label{main thm}
If $(M,H)$ satisfies $\mathrm{(\ast )}$, then
there exists a smooth proper map
$\tilde H\times [0,1]\to G\times\Iso (\mathbb R^d)$,
given by $(\tilde h, t)\to (\psi_t(\tilde h), \e(\tilde h))$,
such that\newline
(i) $(\psi_0, \e)=(\psi, \e)$, \newline
(ii) $\psi_t\co\tilde H\to G$ is a Lie group homomorphism,
for each $t\in [0,1]$, \newline
(iii) $\psi_1(\Pi)$ is finite,\newline
(iv) $(\psi_t,\e )(\Pi)$ is a discrete subgroup of
$G\times\Iso (\mathbb R^d)$ that acts freely
on $\tilde M$,\newline 
(v) $(\psi_t,\e )(\Pi)$ is injective for each $t\in [0,1]$.
\end{thm}
\begin{proof} We first construct $\psi_t\co\tilde H\to G$.
Since $\Pi_1$ is free abelian, the closure of $\psi(\Pi_1)$
is a compact abelian group, which is normalized by $\psi(\tilde H)$
because $\Pi_1\unlhd\tilde H$. Let $T$ be the identity component of
the compact abelian group and let $t$ be the Lie algebra of $T$.
Also $\psi(\tilde H)$ normalizes $T$
so $\Pi_2=\Pi_1\cap\psi^{-1}(T)$ is a normal subgroup of $\tilde H$.
 
Let $\tilde H_Z$ be the centralizer of $\Pi_2$ in $\tilde H$.
Since $\Pi_2\unlhd\tilde H$, we deduce that 
$\tilde H_Z\unlhd\tilde H$. Furthermore, $\tilde H/\tilde H_Z$
is finite. Indeed, $\tilde H_Z$ is the kernel of the $\tilde H$-action 
on $\Pi_2$ be conjugation. 
Since $\Pi_2\subset\tilde H_Z$ and $\Pi_2$ is abelian, 
the $\tilde H$-action 
on $\Pi_2$ be conjugation descends to a Lie group homomorphism
$\tilde H/\Pi_2\to\mathrm{Aut}(\Pi_2)$, where the domain
is compact and the target is discrete, so that the homomorphism
has finite image. 

Let $k=|\tilde H/\tilde H_Z|$.
Let $A=\{\gamma^k\co \gamma\in\Pi_2\}$; this is a normal subgroup
of $\tilde H$.   
Choose a homomorphism $f\co\Pi_2\to t$ with 
$\exp\circ f=\psi_{|\Pi_2}$, and define the map $\phi\co A\to t$
by 
\[\phi(\gamma^k)=\sum_{\tilde h\in\tilde H/\tilde H_Z}
\mathrm{Ad}_{\psi(\tilde h)}(f(\tilde h^{-1}\gamma\tilde h)).\]
It is straightforward to check that $\phi$ is a homomorphism satisfying
$\exp\circ\phi =\psi_{|\Pi_2}$ and 
$\phi(\tilde h a\tilde h^{-1})=\mathrm{Ad}_{\psi(\tilde h)}\phi(a)$ for 
$a\in A$.
Fix an identification of $A$ with the standard lattice 
$\mathbb Z^d\subset\mathbb R^d$,
and use the
$\tilde H$-action on $A$ by conjugation to define the $\tilde H$-action
on $\mathbb R^d$ via
\[\alpha\co\tilde H\to\mathrm{Aut}(\Pi)\cong
\mathrm{GL}_d(\mathbb Z)\subset\mathrm{GL}_d(\mathbb R).\]
Form the semidirect product $\mathbb R^d\rtimes_\a\tilde H$
with multiplication $(v,g)\cdot (w,h)=(v+\a(g)w,gh)$ where 
$g,h\in\tilde H$ and $v,w\in\mathbb R^d$.
We identify $\tilde H$ with the subgroup $\{(0,h)\}$
of $\mathbb R^d\rtimes_\a\tilde H$. 
Define a map
$\Psi\co\mathbb R^d\rtimes_\a\tilde H\to G$ by
\[\Psi\left(\sum_{i}\lambda_i e_i, \tilde h\right)=
\exp\left(\sum_{i} \lambda_i\phi(e_i)\right)\psi(\tilde h),\]
where $\{e_i\}$ is the standard basis in $\mathbb Z^d\subset\mathbb R^d$,
thought of as sitting in $A$, $\lambda_i\in\mathbb R$, and 
$\tilde h\in\tilde H$. Using that $\phi$ is a homomorphism and
$\phi(\tilde h a\tilde h^{-1})=\mathrm{Ad}_{\psi(\tilde h)}\phi(a)$
one checks that $\Psi$ is Lie group homomorphism. Also 
$\Psi_{|\tilde H}=\psi$.

Note that
$N=\{(-a,a)\in \mathbb R^d\rtimes_\a\tilde H\co a\in A\}$
is normal subgroup of $\mathbb R^d\rtimes_\a\tilde H$, and
$N\subset\mathrm{ker}(\Psi)$. Hence $\Psi$ descends
to a Lie group homomorphism $(\mathbb R^d\rtimes_\a\tilde H)/N\to G$ 
Since the subgroups $\mathbb R^d$, $\tilde H$ have trivial intersections
with $N$, they projects isomorphically into 
$(\mathbb R^d\rtimes_\a\tilde H)/N$. Furthermore, $\mathbb R^d$
projects to a normal subgroup whose quotient is isomorphic
to the compact Lie group $\tilde H/A$. Thus, we get an extension
of Lie groups
\[1\to\mathbb R^d\to (\mathbb R^d\rtimes_\a\tilde H)/N\to \tilde H/A\to 1.\]

By~\cite[Proposition 4.3]{Seg} such extensions, with the fixed
$\tilde H/A$-action on $\mathbb R^d$,
are classified by the continuous
cohomology group $H^2_c(\tilde H/A,\mathbb R^d)$. 
Since $\tilde H/A$ is compact, the cohomology group 
vanishes~\cite[section III.2.1]{Gui}, hence
there is exactly one such an extension for each 
$\tilde H/A$-action on $\mathbb R^d$.
Thus the extension coincides with $\mathbb R^d\rtimes_{\bar\a}\tilde H/A$,
where $\bar\a\co \tilde H/A\to\mathrm{GL}_d(\mathbb R)$ 
is the drop of $\a$.
In particular, the composition 
\[I_0\co\tilde H\hookrightarrow \mathbb R^d\rtimes_\a\tilde H\to 
(\mathbb R^d\rtimes_\a\tilde H)/N\cong
\mathbb R^d\rtimes_{\bar\a}\tilde H/A\] 
is an injective Lie group homomorphism.

We write elements of 
$\mathbb R^d\rtimes_{\bar\a}\tilde H/A$
as $(v(\tilde h), \tilde h A)$, where 
$v(\tilde h)\in\mathbb R^d$ and $\tilde h A\in \tilde H/A$. 
Define $I_t\co\tilde H\to\mathbb R^d\rtimes_{\bar\a}\tilde H/A$
by $I_t(\tilde h)=((1-t)v(\tilde h),\tilde h A)$.
This is the deformation of $I_0$ to the homomorphism 
$I_1\co\tilde h\to (0,\tilde h A)$, whose image
is isomorphic to the compact group $\tilde H/A$. 
Note that $I_1(\Pi)=\Pi/A$
is a finite group.

Then $\psi_t(\tilde h)=\Psi\circ I_t$ satisfies (i)--(iii).
To prove (iv) and (v) note that
discreteness of $\Pi_t:=(\psi_t,\e )(\Pi)$
follows from that fact that $\Pi$ is discrete in $\tilde H$
and $(\psi_t,\e )$ is a proper map. Also $\Pi_t$
acts isometrically, and hence properly, on $\tilde M$.
So no infinite order elements of $\Pi_t$ 
can fix a point. If $(\psi_t,\e )(\gamma)$ has finite order for 
$\gamma\in\Pi$, then $(\psi_t,\e )(\gamma)$
is conjugate to $(\psi_0,\e )(\gamma)$ in $G\times\Iso (\mathbb R^d)$
(see e.g.~\cite[Lemma 6.7]{Wil} or~\cite[Lemma 38.1]{CF}), 
so since $\Pi_0$
acts freely, so does $\Pi_t$. Since $\Pi_t$
acts freely, $(\psi_t,\e )$ is injective.
\end{proof}

\begin{rmk}
Similarly to~\cite{Wil}, one can get an explicit bound
on $|\Pi:A|$. We leave it to the interested reader to work out.
Note that our choice of $T$ does not allow to control
$|\Pi_1:\Pi_2|$, however this can be done by
a more delicate choice of $T$ as in~\cite{Wil}. 
\end{rmk}

\begin{thm}\label{thm: deform metric}
If $(M,H)$ satisfies $\mathrm{(\ast )}$, then
there exists a smooth family $q_t$
of complete Riemannian metrics on $M$ such that\newline 
(i) $g_0=g$\newline
(ii) $H\le\Iso(M, g_t)$ for each $t$,\newline 
(iii) there is a smooth family of Riemannian covering maps 
$p_t\co(\tilde M,\tilde g)\to (M, g_t)$,\newline
(iv) there exists a finite index normal subgroup $A$ of 
the group of deck-transformations of $p_1$,
such that the Riemannian covering space 
$(\tilde M/A, \tilde g_1)$ of $(M, g_1)$ 
is isometric to the Riemannian product of $N$, $\mathbb R^{n-d}$, 
and a flat $d$-torus, for some $1\le d\le n$.
\end{thm}
\begin{proof}
Since $\Pi_t\unlhd\tilde H_t$, every 
element $\tilde h_t:=\Psi(I_t((\tilde h))$ of $\tilde H_t$ is 
an $\Pi_t$-equivariant self-diffeomorphism of $\tilde M$, 
hence it descends to 
a self-diffeomorphism $h_t$ of $\tilde M/\Pi_t$.
If $\tilde h_t\in\Pi_t$, then $h_t=\mathrm{id}$, so
the smooth $\tilde H$-action on $\tilde M/\Pi_t$ given by
$\tilde h\to h_t$ factors through
the surjection $\tilde H\to \tilde H/\Pi=H$, and hence descends
to an $H$-action on $\tilde M/\Pi_t$.

We now look at the $\tilde H$-action on $\tilde M\times [0,1]$
given by $\tilde h\cdot(\tilde m, t)=(\tilde h_t(\tilde m), t)$
and the corresponding covering map
$\tilde M\times [0,1]\to (\tilde M\times [0,1])/\Pi$.
The projection $\tilde M\times [0,1]\to [0,1]$ onto the second 
factor descends to the submersion $(\tilde M\times [0,1])/\Pi\to [0,1]$. 
By Morse theory, 
$(\tilde M\times [0,1])/\Pi\to [0,1]$ is a trivial $M$-bundle over $[0,1]$.
By the argument in the previous paragraph, 
the $\tilde H$-action on $\tilde M\times [0,1]$
descends to an $H$-action on the bundle $(\tilde M\times [0,1])/\Pi$
by fiber preserving diffeomorphisms.
Since $H$ is compact, the equivariant covering homotopy 
theorem~\cite{Bie} implies that the bundle
$(\tilde M\times [0,1])/\Pi\to [0,1]$ is $H$-equivariantly smoothly 
isomorphic to the trivial $M$-bundle over $[0,1]$, 
where $H$ acts on $M$ as in $\mathrm{(\ast)}$, and acts trivially 
on $[0,1]$. In particular, for every $t$,
the $\tilde H$-action on $\tilde M/\Pi_t$ given by 
$\tilde h\to h_t$ is smoothly equivalent to the $\tilde H$-action
on $\tilde M/\Pi_0$ given by the surjection $\tilde H\to H$
with kernel $\Pi$.

The composition 
$\tilde M\times [0,1]\to (\tilde M\times [0,1])/\Pi\cong
M\times [0,1]$
is a covering map, which we write as $(\tilde m, t)\to (p_t(m), t)$.
Then $p_t\co\tilde M\to M$ is a covering map with $\Pi_t$
as the group of covering transformations. The map $p_t$ 
is equivariant under the surjection $\tilde H\to H$ given by 
$\tilde h\to \tilde h_t\to h_t\to h_0$, which 
is simply the surjection with kernel $\Pi$.

Since $\Pi_t$ acts isometrically on $\tilde M$, 
there is a unique Riemannian metric $g_t$ on $M$ 
that makes $p_t$ a Riemannian covering. 
Thus $p_t$ is a smooth
family of Riemannian coverings and $g_t$ is a smooth
family of Riemannian metrics. Also $H\le\Iso(M,g_t)$,
and since the group $\psi_1(A)$ is trivial, $(M, q_1)$
has a finite Riemannian cover isometric to $\tilde M/A$,
which is the Riemannian product of $N$, $\mathbb R^{n-d}$, 
and a flat $d$-torus.
\end{proof}
\begin{rmk}
\label{h-cobordism} 
Let $M$ be a compact manifold with infinite fundamental group
that admits a metric of $\Ric\ge 0$, so that 
$\pi_1(M)$ is virtually-$\mathbb Z^n$ and
the universal cover of $M$ is isometric to 
$C\times\mathbb R^n$~\cite{CG2} with $\Ric(C)\ge 0$. 
The choice of a metric of $\Ric\ge 0$ on $M$ uniquely specifies the 
isometry type of $C$, and hence the diffeomorphism type of $C$.
It is natural to ask whether $M$ could admit a different metric,
with the universal cover 
of $M$ isometric to $C^\prime\times\mathbb R^n$, 
for which $C^\prime$ and $C$ are non-diffeomorphic.

It turns out that the $s$-cobordism theorem rules out this
possibility if $\dim(C)\ge 5$. Namely, a finite cover of $M$
has to be diffeomorphic to $C\times T^n$ as well as
to $C^\prime\times T^n$. 
The diffeomorphism $\phi\co C\times T^n\to
C^\prime\times T^n$, defines an automorphism of $\pi_1(T^n)$, and hence
a self-homotopy equivalence of $T^n$.
Since any self-homotopy equivalence of $T^n$ is
homotopic to a diffeomorphism, we can precompose $\phi$
by this diffeomorphism, so that we can assume that
$\phi$ induces the trivial homomorphism of the fundamental groups.
If $n=1$, then the lift of $\phi$ to the universal covers
defines an $h$-cobordism between $C$ and $C^\prime$, so they are
diffeomorphic. If $n>1$, we can proceed by induction
splitting off one circle at a time, and getting an $h$-cobordism
between $C\times T^k$ and $C^\prime\times T^k$, which is trivial
because $\mathbb Z^k$ has trivial Whitehead torsion.
After splitting all the circle factors, we again conclude that 
$C$ and $C^\prime$ are diffeomorphic.

Note that there do exist simply-connected $4$-manifolds
$C$, $C^\prime$ that become diffeomorphic after taking products
with $S^1$. For example, by Seiberg-Witten theory,
the $K3$-surface has infinitely many
smooth structures, but each compact manifold of dimension $\ge 5$
has only finitely many smooth structures. It it unclear however
if such examples can admit $\Ric\ge 0$.
\end{rmk}

\section{Equivariant connected sum and positive scalar curvature}
\label{sec: pos scal curv}

Let $M$, $N$ be manifolds of the same dimension
each equipped with a smooth action of
a compact Lie group $G$ that has a nonempty set of the fixed points.
Assume that at some fixed points $m\in M$, $n\in N$
the isotropy representations of $G$ are equivalent.
If the $G$-actions are isometric in some Riemannian
metrics on $M$, $N$, then there are small balls around $m$, $n$
which are $G$-equivariantly diffeomorphic. The {\it equivariant connected
sum} of $M$ and $N$ at the points $m$, $n$ is a $G$-manifold
obtained by removing the above balls and gluing their
complements by a $G$-equivariant diffeomorphism of the boundary 
spheres. For the purposes of this paper we ignore the issue
of orientability for $M\# N$, so we need not put any restrictions on the
diffeomorphism of the boundary spheres.

It follows easily from~\cite{GL} that 
positivity of scalar curvature is preserved under equivariant 
connected sums. Namely, let $(M,g)$ be a complete manifold of $\scal>0$
of dimension $m\ge 3$,
and let $G$ be a compact subgroup of $\Iso(M)$ that fixes a point
$x\in M$. It was shown in~\cite{GL} that $M\setminus\{x\}$ has a 
metric $g^\prime$ of positive scalar curvature for which 
there exist small $R>r>0$ such that
outside the ball $B_g(x,R)$ we have $g^\prime=g$, and 
on $B_g(x,r)\setminus\{x\}$ the metric
$g^\prime$ is the product $\mathbb R\times S^{m-1}$, 
where $S^{m-1}$ is round sphere of a small fixed radius. 
Furthermore, the metric $g^\prime$
on $B_g(x,R)\setminus\{x\}$ is rotationally symmetric, in particular,
the $G$-action on $M\setminus\{x\}$ is isometric 
with respect to $g^\prime$. Thus, 
one can keep the scalar curvature positive while taking
equivariant connected sums.

In our applications we take equivariant connected sums of a fixed 
$G$-action on $M$ and a certain $G$-action on the standard sphere 
$S^m$, whose set of the fixed points is a homology sphere $S$. 
Thus, the equivariant connected sum $M\# S^m$ is diffeomorphic to $M$,
while the $G$-action changes considerably.
Namely, if $F$ is the component of the fixed point set of the
original $G$-action on $M$ that
contains the point at which we take the equivariant connected sum, 
then the corresponding component of the fixed point set 
of the $G$-action on $M\# S^m$ is diffeomorphic to $F\# S$.

First, we need background on homology spheres that bound
(smooth) contractible manifolds. In each dimension 
$n\ge 3$ there are non-simply-connected 
homology $n$-spheres that bound contractible manifolds.
In fact if $n\ge 5$, then 
after changing a smooth structure any homology sphere
bounds a contractible manifold~\cite{Ker}. Also any
homology $4$-sphere bounds a contractible manifold~\cite{Ker}.
Some homology $3$-spheres bound contractible manifolds~\cite{Gla, Gor} 
while some do not~\cite{Ker}. (Note that any homology $3$-sphere
bounds a topological contractible manifold~\cite{Fre}, but
we need to work in smooth category). 
In our application, $\partial(D^k\times C)$ must be 
diffeomorphic to $S^{k+n}$, where
$D^k$ denotes the closed $k$-disk with $k>0$. 
It is easy to see that $\partial(D^k\times C)$ is simply-connected,
so if $\dim(D^k\times C)\ge 6$ then by the $s$-cobordism
theorem applied to $D^k\times C$ with a small open ball in 
$\mathrm{Int}(D^k\times C)$
removed, $\partial(D^k\times C)$ is diffeomorphic to $S^{k+n}$.
If $D^k\times C$ is a $5$-manifold, so that $k=1$, $\dim(C)=4$,
the above argument does not work, because the smooth Poincar\'e
conjecture is still open. So we have to
restrict ourselves only to those $C$'s 
for which $\partial(D^k\times C)$ is diffeomorphic to $S^4$. 
This is true in examples of~\cite{Gla} 
(note that $\partial(D^k\times C)$ is the double of $C$).

Now fix a homology $n$-sphere $S$
that bounds a (smooth) contractible manifold $C$,
and if $n=3$ we assume that the double of $C$ is diffeomorphic to $S^4$.
Let $G$ be a compact Lie group acting on $\mathbb R^k$ via
an embedding $\rho\co G\hookrightarrow O(k)$, and suppose that
the only fixed point of $\rho(G)$ is the origin. 
Consider the closed unit ball $D^k$ in $\mathbb R^k$ centered at the
origin with the induced $G$-action. Also assume that $G$ acts
trivially on $C$.
Given the data, we next construct a $G$-action on $S^m$, with $m=n+k$, 
such that the fixed point set is $S$ and the isotropy representation
at any fixed point is equivalent to the direct sum of $\rho$
and the trivial $n$-dimensional representation.

Look at the product of the $G$-actions on 
on $D^k\times C$. By equivariantly smoothing the corners 
we get a smooth $G$-action on the manifold
$\partial(D^k\times C)=(D^k\times S)\cup (\partial D^k\times C)$ 
which is diffeomorphic to $S^m$.
Since the origin is the only fixed point of $\rho(G)$,
the set of the fixed points for the $G$-action on $\partial(D^k\times C)$
is the homology sphere $\{0\}\times S$, and 
the isotropy representation at any
fixed point is equivalent to the direct sum of $\rho$ and
the trivial $n$-dimensional representation.

By taking equivariant connected sums of $q$ copies of the 
$G$-action on $S^m$, we get the $G$-action on $S^m$ with the same
isotropy representation at a fixed point and such that
the set of the fixed points is the connected sum $\#_q S$
of $q$ copies of $S$. These actions are topologically
inequivalent for $\#_q S$ and $\#_p S$ are not homeomorphic 
if $q\neq p$, because  by Grushko's theorem~\cite{Sta}
the rank (i.e. the minimum number of generators) of 
$\pi_1(\#_q S)\cong\ast_q\pi_1(S)$ is 
$q\cdot\mathrm{rk}(\pi_1(S))$.

Finally, we show that if $n\ge 5$, then 
these actions on $S^m$
can be arranged to preserve some metrics of $\scal>0$.
By Lemmas below we can assume that $C$ has a metric of $\scal>0$
that is the product metric near the boundary. Equip $D^k$ with an
$O(k)$-invariant metric of $\scal\ge 0$ which is the product
metric near the boundary. Then the product
metric on $D^k\times C$ has $\scal>0$ away from the boundary.
By smoothing the metric at the corners as in~\cite{Gaj}, we
can assume that $\partial(D^k\times C)$ has a $G$-invariant
metric of $\scal>0$.

\begin{lem} Let $\pi$ be a finitely presented group with
$H_1(\pi)=0=H_2(\pi)$ and let $n\ge 5$ be an integer.
Then there exists a homology $n$-sphere with fundamental 
group $\pi$ and of $\scal>0$. 
\end{lem}
\begin{proof}
By~\cite[p.68]{Ker} the desired homology $n$-sphere can be constructed from
a connected sum of several copies of $S^1\times S^{n-1}$ by surgeries
on circles and $2$-spheres. 
Thus the surgeries have codimension 
$\ge 3$ and the homology sphere have $\scal>0$ by~\cite{SY, GL}.
\end{proof}

\begin{lem} If $n\ge 5$, then
after possibly changing a smooth structure,
each homology $n$-sphere $M$ of $\scal(M)>0$
bounds a contractible manifold with $\scal>0$
so that the metric is a product near the boundary.
\end{lem}
\begin{proof}
We argue as in~\cite[pp.70--71]{Ker}. First by doing 
surgery on circles and $2$-spheres we turn $M$ into a homotopy sphere $\Sigma$
of $\scal>0$. By~\cite{Gaj} the trace of this surgery (a cobordism $W$
between $M$ and $\Sigma$)
also has $\scal>0$ so that the metric is the product near the boundary. 
We think of $W$ as $M\times I$ with some handles attached to
$M\times\{1\}$, where $I=[0,1]$. For some $m\in M$, the interval
$m\times I$ touches none of the previous surgeries and 
meets orthogonally
each copy of $\partial W\times\{t\}$ near the boundary, where
the product in $\partial W\times\{t\}$ refers to the product metric.
Let $N(m\times I)$ be a small open tubular neighborhood of $m\times I$.
Then following~\cite{Gaj}, one can smooth the metric at the corners of 
$W^\prime=W\setminus N(m\times I)$ to a
metric that has $\scal>0$ and is a product near the
boundary. As in~\cite[p.71]{Ker}, $W^\prime$ is a contractible manifold
with boundary $M\#\Sigma$.  
\end{proof}

\section{Constrains on isometric actions}

Part (iv) of Theorem~\ref{thm: deform metric} yields a nontrivial 
restriction on the topology $M$, namely, a finite cover of $M$ splits
off a torus factor. This idea was explored in~\cite{BK1, BK2}, 
to show, for example, that many vector bundles admit no complete metric of 
$\sec\ge 0$. 

The purpose of this section is to show that 
Theorem~\ref{thm: deform metric} also 
gives a lot of information information on compact subgroups of $\Iso(M)$.
This leads, in particular, to numerous examples of group 
actions on compact manifolds that are non-isometric in any 
metric of $\Ric\ge 0$, as well as examples of group 
actions on manifolds that are non-isometric in any 
complete metric of $\sec\ge 0$.

Our first goal is to lift the the isometric $H$-action on 
$(M,g_1)$, constructed in Theorem~\ref{thm: deform metric},
to the covering $\tilde M/A\to M$. This is generally impossible, 
e.g. the circle action on itself by left translations does
not lift to any cover. 
However, by the standard results on lifting groups actions to 
covering spaces, which are reviewed in 
Appendix~\ref{lift group act}, the $H$-action on $M$ always 
lifts to an $\tilde H/A$-action on $\tilde M/A$,
and furthermore, the following is true.
\begin{enumerate}
\item[(1)] If $H$ fixes a point $x\in M$, then the $H$-action on $M$
lifts to an $H$-action on $\tilde M/A$.
\item[(2)]
If $H$ is connected, the  $H$-action on $M$ lifts to an 
$\bar H$-action on $\tilde M/A$, where $\bar H$ is
a connected Lie group that surjects onto $H$ and the kernel
of the surjection is a subgroup of $\Pi/A$, 
which is a finite group.
\item[(3)] 
If $H$ is connected and simply-connected, 
then any $H$-action on $M$
lifts to an $H$-action on $\tilde M/A$.
\end{enumerate}
Since $(\tilde M/A, \tilde g_1)$ is isometric to 
$N\times\mathbb R^{n-d}\times T^d$,
where $T^d$ is a flat torus, 
in the cases (1), (3), $H$ acts on $\tilde M/A$
via an embedding 
into $\Iso(N)\times\Iso(\mathbb R^{n-d})\times \Iso(T^d)$.
In particular, the $H$-action 
preserves the decomposition of $\tilde M/A$ into the product of
$N$, $\mathbb R^{n-d}$, and $T^d$.
Of course, isometric actions of compact Lie groups on 
$\mathbb R^{n-d}$ and $T^{d}$ are well-understood. 
Using the results on actions on tori, stated in 
Section~\ref{actions on tori}, we see that 
in the cases (1), (3), the $H$-action 
on the $T^{d}$-factor is as follows.
\begin{itemize}
\item[$\mathrm{(1)^\prime}$]
Assume $H$ fixes a point $x\in M$, so the lifted
$H$-action on $\tilde M/A$ also fixes a point, and let $t$ be
the fixed point of the $H$-action on the $T^{d}$-factor.
Then if $h\in H$ acts nontrivially on the $T^{d}$-factor, then the
induced action of $h$ on $\pi_1(T^{d},t)$ is nontrivial. 
In particular, if $H$ acts trivially on $\pi_1(M,x)$, then 
the $H$-action on the $T^{d}$-factor is trivial.
\item[$\mathrm{(3)^\prime}$]
If $H$ is connected and simply-connected,
then the $H$-action on the $T^{d}$-factor is trivial.
\end{itemize}

\begin{rmk} \label{rmk: t^d factor}
The fact that $H$ acts trivially
on the $T^d$-factor implies in particular that 
each component of the fixed points set
has a direct $T^d$-factor. This is the key 
point in Theorem~\ref{eq-sum-star} below.
\end{rmk}

\begin{thm}\label{eq-sum-star} Let $H$ be a compact
Lie group acting smoothly and effectively on a manifold $M$.
Suppose the set of the fixed 
points of $H$ has a component $F$ of dimension $s\ge 3$ such that
for some $x\in F$, the $H$-action on $\pi_1(M, x)$ is trivial.
Then there are infinitely many smooth $H$-actions on $M$
for which $(M, H)$ does not satisfy $\mathrm{(\ast)}$.
\end{thm}
\begin{proof}
Following Section~\ref{sec: pos scal curv}, we construct a smooth $H$-action
on the standard sphere $S^{m}$, where $m=\dim(M)$, 
such that the fixed point set is a non-simply-connected 
integral homology $s$-dimensional sphere $S$, and the isotropy 
actions at $x\in M$ and at some fixed point of $S^m$ are equivalent.
By replacing $(S^m, H)$ with the equivariant connected sum at $x$
of the two copies of $(S^m, H)$, we can assume that $S$ is a connected
sum of two non-simply-connected homology spheres. In particular, 
$\pi_1(S)$ has infinitely many ends~\cite[4.A.6.6]{Sta}.
Then the equivariant connected sum of $M$ and  $S^{m}$ at the
fixed points is a smooth $H$-action on $M$ such that
one of the components $F^\prime$ of the fixed point  set 
is diffeomorphic $F\# S$. Also the induced $H$-action on 
$\pi_1(M,y)$ is trivial for $y\in F^\prime$. 

Arguing by contradiction, assume 
the new $H$-action on $M$ satisfies $\mathrm{(\ast)}$,
for some $\tilde M$, $\Pi$, etc. 
Now apply Theorem~\ref{main thm} to find $A$ and $g_t$, and
let $q\co \tilde M/A\to M$ be the corresponding covering.
Pick $\tilde x\in \tilde M/A$ with $q(\tilde x)=x$, and let
$\bar F$ be the component of $q^{-1}(F^\prime)$ that contains $\bar x$.
The $H$-action on $M$ lifts uniquely to an $H$-action on $\tilde M/A$
that fixes $\bar x$. Then the set of the fixed points of the $H$-action
on $\tilde M/A$ has $\bar F$ as a component.
Note that
$q_{|\bar F}\co \bar F\to F^\prime$ is a finitely-sheeted 
covering. By $\mathrm{(3)^\prime}$ above, 
$\bar F$ is the product of $T^{d}$ and a submanifold of $N$,
and therefore, $\pi_1(\bar F)$ 
has a direct $\mathbb Z^d$-factor.

Hence the group $\pi_1(F^\prime)$ has $1$ or $2$ 
ends~\cite[4.A.6.1, 4.A.6.3]{Sta}, where $\pi_1(F^\prime)$ has $2$ ends
precisely when $d=1$ and $\pi_1(\tilde F)/\mathbb Z^d$ is finite. 
On the other hand, $\pi_1(F^\prime)=\pi_1(F)\ast\pi_1(S)$ 
so $\pi_1(F^\prime)$ has infinitely many ends~\cite[4.A.6.6]{Sta}.

To get infinitely many examples of such $H$-actions we proceed 
as follows. Fix a homology $s$-sphere $S$ that bounds a contractible 
manifold, whose double is $S^4$ if $s=3$. 
For each $s$-dimensional connected component $F$ 
of the fixed point set of $(M,H)$, we build an $H$-action on $S^m$
with fixed point set $S$ and the 
isotropy representation at a fixed point
equivalent to the isotropy representation at a point $f$ of $F$.
For every positive integer $q$ and every $F$, we 
take the equivariant connected sum of $M$ at $f$ with $q$ copies 
of $(S^m,H)$, so that for this new $H$-action $\rho_q$ on $M$,
each $F$ turns into $F_q:=F\# q S$. 
By Grushko's theorem~\cite{Sta}, 
the rank (i.e. the minimum number of generators)
of $\pi_1(F_q)$ is $\mathrm{rk}(\pi_1(F))+q\cdot\mathrm{rk}(\pi_1(S))$.
Let $r_q$ be the minimum of the numbers $\mathrm{rk}(\pi_1(F_q))$
over all $s$-dimensional
components of the fixed point set of $\rho_q$.
The minimum is of course realized on some component, so
if $q\neq p$, then $r_q\neq r_p$, therefore, 
$\rho_q$ and $\rho_p$ are topologically inequivalent.
\end{proof}
\begin{proof}[Proof of Theorem~\ref{ric-scal-eq-sum}]
We construct the $H$-actions on $M$ as in the proof of
Theorem~\ref{eq-sum-star}, in particular,
the actions preserve no metrics of $Ric\ge 0$.
If $\scal(M)>0$ and $\dim(F)\ge 5$, 
then by Section~\ref{sec: pos scal curv},
the $H$-action obtained by the equivariant connected sum as above
preserves a metric of $\scal>0$.
\end{proof}
\begin{proof}[Proof of Theorem~\ref{sec-scal-eq-sum}]
The argument here is very similar to
the proof of Theorem~\ref{eq-sum-star}, in fact we only 
need to replace the second paragraph in that proof.
Namely, we need another argument explaining that 
some finite index subgroup of $\pi_1(F^\prime)$ has a 
$\mathbb Z^d$-factor for $d\ge 1$.
 
Since $S$ is a nontrivial connected sum,
$\pi_1(S)$ is infinite which ensures that 
$\pi_1(F^\prime)\cong\pi_1(F)\ast\pi_1(S)$ is infinite.
Arguing by contradiction assume that the new $H$-action
preserves a complete metric of $\sec\ge 0$. Then $F^\prime$
is a totally geodesic submanifold, and hence $\sec(F^\prime)\ge 0$.
Therefore, a finite index subgroup of $\pi_1(F^\prime)$ has a 
$\mathbb Z^d$-factor for $d\ge 1$. 

If $\scal(M)>0$ and $\dim(F)\ge 5$, 
then by Section~\ref{sec: pos scal curv},
the $H$-action obtained by the equivariant connected sum as above
preserves a metric of $\scal>0$.
\end{proof}

\begin{rmk}
Theorem~\ref{sec-scal-eq-sum} can be extended 
to complete manifolds of
nonnegative $k$-Ricci curvature (see~\cite{Wu} for a definition),
provided $\dim(F)\ge k+1$, because these assumptions force $F$ 
to have nonnegative Ricci curvature, so that $\pi_1(F)$ is 
virtually free abelian. The only $k$-Ricci curvature 
for which the extension is void is the Ricci curvature
itself, i.e. the $k$-Ricci curvature with $k=\dim(M)-1$;
indeed, $k=\dim(M)-1$ implies $M=F$. In general, little is known
about the totally geodesic submanifolds in manifolds of 
$\Ric\ge 0$, which is one of the obstacles to understanding
simply-connected manifolds with $\Ric\ge 0$.  
\end{rmk}

\appendix

\section{Lifting group actions to covering spaces.}
\label{lift group act}
If $p\co\bar X\to X$
is a covering and $\pi\co\bar G\to G$ is a surjection of Lie groups,
then the $\bar G$-action $\bar\theta\co\bar G\times \bar X\to\bar X$
on $\bar X$ is said to {\it cover} the $G$-action $\theta\co G\times X\to X$
on $X$ if $p\circ\bar\theta=\theta\circ (\pi\times p)$.
We also say that $\theta$ {\it lifts} to $\bar\theta$.

It is shown in~\cite[section I.9]{Bre} that if $p\co\tilde X\to X$
is a regular covering with deck transformation group $\Delta$, 
and if the Lie group $G$ acts effectively on $X$, then the set of all lifts
of elements of $G$ to self-homeomorphisms of $\tilde X$ is 
a Lie group $\tilde G$ that surjects onto $G$ with kernel 
$\Delta$, so that the $\tilde G$-action on $\tilde X$ covers
the $G$-action on $X$. 
In fact, in many cases one can find a much smaller Lie subgroup
$\bar G$ of $\tilde G$ such that 
the $\bar G$-action on $\tilde X$ still covers
the $G$-action on $X$. The kernel of the surjection $\bar G\to G$
is a subgroup of $\Delta$.

For example, if $G$ fixes a point
$x\in X$ and the induced $G$-action on $\pi_1(X,x)$ 
preserves the subgroup $p_*(\pi_1(\tilde X,\tilde x))$
for $\tilde x\in p^{-1}(x)$,
then we can take $\bar G=G$, because then, 
any $g\in G$ lifts uniquely to $\tilde g\in \tilde G$ with
$\tilde g(\tilde x)=\tilde x$~\cite[section I.9]{Bre}.

If $G$ is connected, then one can choose $\bar G$ to be 
connected~\cite[section I.9]{Bre}, so that $\bar G$ is a covering
space of $G$.  In particular, if $G$ is also
simply-connected, then we can take $\bar G=G$, in other words,
any action of a connected simply-connected Lie group $G$ lifts
to a $G$ action on $\tilde X$. Also if $G$ is torus and 
$p\co\tilde X\to X$ is a finite cover, then $\bar G$ is also
a torus and the surjection $\bar G\to G$ is a finite cover.

\section{Actions on tori.}
\label{actions on tori}
If a connected compact Lie group $K$ acts effectively on a torus $T^m$, 
then $K$ is a torus, the action is free, and 
in fact $T^m$ is a trivial principal $K$-bundle~\cite[section IV.9]{Bre}.
It seems the base $T^m/K$ of the bundle must also be diffeomorphic to
a torus, but I do know how to prove it. 
The result is clear if the $K$-action preserves a flat metric,
because in this case, $T^m/K$ carries a Riemannian submersion metric
of $\sec\ge 0$, and by the homotopy sequence of the fibration
$T^m/K$ has contractible universal cover and by above 
$\pi_1(T^m/K)$ is free abelian, so the metric on 
$T^m/K$ is flat~\cite{CG1}, and hence $T^m/K$ is a torus~\cite{Wol}.

If a compact Lie group $K$ acting effectively
on a torus $T^m$ fixes a point $t$,
then the induced $K$-action on $\pi_1(T^m,t)$ is 
effective~\cite[section IV.9]{Bre}.
Here is an example of an involution on $T^m$ that is not isometric in
any metric of $\scal\ge 0$.
Consider an isometric involution of
$T^m$ whose fixed-point-set has a component $F$ of dimension $\ge 3$. 
Then by taking 
an equivariant connected sum at a point of $F$ as in 
Section~\ref{sec: pos scal curv}, we get an involution
on $T^m$ such that a component
of the fixed point set is diffeomorphic to $F\# S$ for
a homology sphere $S$.
Now if the action is isometric in a flat metric, then
$F\# S$ is totally geodesic, hence flat, so $\pi_1(F\# S)$
is virtually abelian, but $\pi_1(F\# S)\cong\pi_1(F)\ast\pi_1(S)$ 
is not virtually abelian
unless $S$ is simply-connected. Since any metric of $\scal\ge 0$
on a torus is flat~\cite{LM}, 
we have constructed an involution on a torus
that is not isometric in any metric of $\scal\ge 0$.

\small
\bibliographystyle{amsalpha}
\bibliography{eq}

\

DEPARTMENT OF MATHEMATICS, 253-37, CALIFORNIA INSTITUTE OF TECHNOLOGY,
PASADENA, CA 91125, USA

{\normalsize
{\it email:} \texttt{ibeleg@its.caltech.edu}}

\end{document}